\documentclass[british,english]{article}
\usepackage[T1]{fontenc}
\usepackage[latin9]{inputenc}
\usepackage{color}
\usepackage{amsmath}
\usepackage{amssymb}
\usepackage{babel}

\newtheorem{theo}{Theorem}

\newtheorem{coro}{Corollary}

\newtheorem{rema}{Remark}

\def\be#1\ee{\begin{equation}#1\end{equation}}


 \global\long\def\sbr#1{\left[ #1\right] }
 \global\long\def\cbr#1{\left\{  #1\right\}  }
 \global\long\def\rbr#1{\left(#1\right)}

 \global\long\def\R{\mathbb{R}}
 
 \global\long\def\N{\mathbb{N}}
 \global\long\def\Z{\mathbb{Z}}

 \global\long\def\dd#1{\textnormal{d}#1}

 \global\long\def\ra{\rightarrow}
 
 \global\long\def\ns{\infty}

\begin{document}

\selectlanguage{english}%

\title{The It{\^o}-F\"ollmer formula -- nonstandard cases}

\author{W. M. Bednorz, R. M.  \L ochowski, P. L. Zondi, F. J. Mhlanga and D. Hove}

\maketitle

\begin{abstract}
The purpose of this note is to prove the It{\^o}-F\"ollmer formula for
the c\`adl\`ag paths possessing quadratic variation in a possibly ``weakest''
sense along some sequence of partitions. By this we mean, for example,
that we do not require the jumps of the quadratic variation to occur exactly
at the times when the path itself jumps or that the quadratic
variation exists at all time instances. Moreover, we deal with more
general partitions than partitions with vanishing mesh used by F\"ollmer, also relaxing restrictions imposed on a sequence of partitions in other related literature.  
\end{abstract}

\section{Introduction, notation and definitions}

By $\Z$ we will denote the set of all integers, by $\N$ the set
of all positive integers and by $\N_{0}$ the set of all non-negative
integers. 

Let $T$ be a positive real number. By a \emph{partition} of the interval
$(0,T]$ we will mean a finite sequence $\pi=\rbr{(t_{i},t_{i+1}]}_{i=0}^{k-1}$
of intervals $(t_{0},t_{1}]$, $(t_{1},t_{2}]$, ..., $(t_{k-1},t_{k}]$,
such that 
\[
0=t_{0}<t_{1}<\ldots<t_{k}=T.
\]
A sequence of partitions (of the interval $(0,T]$) denoted by $\rbr{\pi^{n}}$ or
$\rbr{\pi(n)}$ will be a sequence such that for any $n\in\N$, $\pi^{n}$
(or $\pi(n)$) is a partition of the interval $(0,T]$. 

The celebrated It{\^o}-F\"ollmer formula proved in \cite{Foellmer:1981}
assumes that a deterministic c\`adl\`ag path $x:[0,T]\ra\R$ possesses
quadratic variation along some sequence $\rbr{\pi^{n}}=\rbr{(t_{i-1}^{(n)},t_{i}^{(n)}]_{i=1}^{k(n)}}$
of partitions of the interval $(0,T]$ in the sense that
\begin{enumerate}
\item the measures
\[
\mu_{n}:=\sum_{i=1}^{k(n)}\rbr{x\rbr{t_{i}^{(n)}}-x\rbr{t_{i-1}^{(n)}}}^{2}\delta_{0}\rbr{\cdot-t_{i-1}^{(n)}}=\sum_{(u,v]\in\pi^{n}}\rbr{x\rbr v-x\rbr u}^{2}\delta_{0}\rbr{\cdot-u}
\]
tend weakly to some non-negative, finite, Borel measure $\mu$ on
$[0,T]$;
\item the atoms of the measure $\mu$ are exactly the sets consisting of
times where the jumps of the path $x$ occur, and have the same mass
as the squares of the jumps of $x$ at these times: for any $t\in(0,T]$
\[
\mu\rbr{\cbr t}=\rbr{\Delta x(t)}^{2},
\]
where $\Delta x(t):=x(t)-x(t-).$
\end{enumerate}
By the Portmanteau Theorem, see \cite[Theorem 2.1]{Billingsley:1999cl}, Assumption 1. is equivalent to the fact that the pointwise limit of the cumulative distribution functions
(CDFs) of the measures $\mu_{n}$, $n\in\N$, exists at any point
$t\in[0,T]$ where the CDF of $\mu$ is continuous and at the point
$t=T$, and is equal $\mu[0,t]$:
\begin{equation}
\lim_{n\ra+\ns}\mu_{n}[0,t]=\lim_{n\ra+\ns}\sum_{(u,v]\in\pi^{n},\ u\le t}\rbr{x\rbr v-x\rbr u}^{2}=\mu[0,t].\label{eq:qvar}
\end{equation}
Assumption 2. yields that the CDF of $\mu$ is continuous at any point $t$
where $x$ is continuous. F\"ollmer also assumed that 
\begin{equation}
\text{mesh}\rbr{\pi^{n}}:=\max_{1\le i\le k(n)}\rbr{t_{i}^{(n)}-t_{i-1}^{(n)}}\ra0\text{ as }n\ra+\ns\label{eq:mesh}
\end{equation}
and as a result of this, the c\`adl\`ag property of $x$ and the monotonicity
of $t\mapsto\mu_{n}[0,t]$, the pointwise limit of the CDFs of the measures $\mu_{n}$ exists at
\emph{any} point $t\in[0,T]$ and is equal $\mu[0,t]$. 

From (\ref{eq:mesh}) it follows that $\mu[0,t]$ is also equal to
the limit 
\begin{equation}
\lim_{n\ra+\ns}\sum_{(u,v]\in\pi^{n}}\rbr{x\rbr{v\wedge t}-x\rbr{u\wedge t}}^{2}\label{eq:qvar1}
\end{equation}
(see for example \cite[Sect. 2.1]{LochOblPS:2021}). The limit (\ref{eq:qvar1}), whenever
it exists, will be called the \emph{quadratic variation of $x$ till
time $t$ along the sequence of partitions $\rbr{\pi^{n}}$}.

In the classical stochastic calculus and in F\"ollmer's setting the
jumps of the quadratic variation occur at the points where $x$ has
jumps only. However, it is well known that there are counterexamples
of continuous functions $x$ such that (\ref{eq:qvar1}) is not a
continuous function of $t$, see \cite[Sect. 5.2]{Coquet2006}, it may be actually any non-decreasing
function of $t$ starting from $0$; see \cite[Theorem 7.1]{Obloj_local:2015}.

In this note we investigate under what assumptions the convergence
of the sums 
\[
\sum_{(u,v]\in\pi^{n}}f''\rbr{x\rbr u}\rbr{x\rbr v-x\rbr u}^{2}\text{ and }\sum_{(u,v]\in\pi^{n}}f'\rbr{x\rbr u}\rbr{x\rbr v-x\rbr u}
\]
holds for any twice continuously differentiable $f:\R\ra\R$. We also
investigate in which cases the limit (assuming that it exists) 
\[
\lim_{n\ra+\ns}\sum_{(u,v]\in\pi^{n}}f''\rbr{x\rbr u}\rbr{x\rbr v-x\rbr u}^{2}
\]
may be represented as an integral with respect to some finite Borel
measure $\mu$ on $[0,T]$ and a version of the It{\^o}-F\"ollmer formula
holds. We will not assume that the atoms of this measure consist
of the points of jumps of $x$ only. We will work with more general sequences
of partitions $\rbr{\pi^{n}}$ than those considered by F\"ollmer, since the partitions
with vanishing mesh do not contain for example the Lebesgue partitions
of paths which are constant on some interval. 

A definition of the \emph{Lebesgue partition} $\pi^{c,r}$ (sometimes,
for typographical reasons, we will also denote it by $\pi(c,r)$)
of the interval $[0,T]$ for continuous $x:[0,T]\ra\R$ and the grid
\[
c\cdot\Z+r=\cbr{z\in\R:\exists p\in\Z,z=p\cdot c+r},
\]
where $c>0$ and $r\in[0,c)$ is the following. First we define $\tau_{0}^{c,r}=0$,
and for $k\in\N$, 
\[
\tau_{k+1}^{c,r}:=\begin{cases}
\rbr{\inf\cbr{t>\tau_{k}^{c,r}:x_{t}\in\rbr{c\cdot\Z+r}\setminus\cbr{x_{\tau_{k}^{c,r}}}}}\wedge T & \text{if }\tau_{k}^{c,r}<T,\\
T & \text{if }\tau_{k}^{c,r}=T,
\end{cases}
\]
(we apply the convention that $\inf\emptyset=+\ns$). Denote $k(c,r)=\min\cbr{k\in\N:\tau_{k}^{c,r}=T}$.
Now we define $\pi^{c,r}$ as the family of intervals 
\[
\pi^{c,r}=\cbr{\left(\tau_{k}^{c,r},\tau_{k+1}^{c,r}\right],k=0,1,\ldots,k(c,r)-1}.
\]
A definition of the Lebesgue partitions for general c\`adl\`ag paths may
be found for example in \cite[Sect. 4]{Vovk_cadlag:2015}.

Our toy example will be a similar zigzag function to the one from
\cite[Sect. 5.2]{Coquet2006}. We will denote it by $z$. $z$ is defined as follows.
$z:[0,1]\ra\R$, $z(t)$ is equal to $0$ at each $t=1-2^{-2m}$,
equal $(m+1)^{-1/2}$ at each $t=1-2^{-2m-1}$, $m\in\N_{0}$, and
is affine between these points. However, instead of dyadic partitions
we will consider the Lebesgue partitions. For $\alpha\in[0,1)$ let
$\rbr{\rho^{n}(\alpha)}$ denotes the sequence of the Lebesgue partitions
$\rbr{\pi\rbr{1/\sqrt{n},\alpha/\sqrt{n}}}$ for $z$. It is possible
to prove that $\text{mesh}\rbr{\rho^{n}(\alpha)}\ra0$ as $n\ra+\ns$ and
$y$ has discontinuous quadratic variation along the sequence $\rbr{\rho^{n}(\alpha)}$
for any $\alpha\in[0,1)$. More precisely, it is possible to prove
that 
\begin{equation} \label{zigzag}
\lim_{n\ra+\ns}\sum_{(u,v]\in\rho^{n}(\alpha)}\rbr{z\rbr{v\wedge t}-z\rbr{u\wedge t}}^{2}=
\begin{cases}
0 & \text{if } \, t \in [0, 1), \\
L(\alpha) & \text{if } \, t =1,
\end{cases}
\end{equation}
where
\[
L(\alpha):=\lim_{n\ra+\ns}\frac{2}{n}\sum_{m=1}^{n}\left\lfloor \sqrt{\frac{n}{m}}-\alpha\right\rfloor =2\sum_{l=1}^{+\ns}\cbr{\frac{l}{(\alpha+l)^{2}}-\frac{l}{(\alpha+l+1)^{2}}}
\]
($\lfloor a\rfloor$ is the floor of $a$), see \cite{PhumlaniPhD} or the Appendix. $L$ is continuous and decreasing. $L(0)=\pi^{2}/3\approx3.29$, $L(1/2)=\pi^{2}-8\approx1.87$
while $\lim_{\alpha\ra1-}L(\alpha)=\pi^{2}/3-2\approx1.29$. 

Let us continue with further definitions. In what follows we assume
that $x:[0,T]\ra\R$ is c\`adl\`ag. 

By $D$ we denote the set of points where $x$ has jumps 
\[
D=\cbr{t\in(0,T]:\left|\Delta x(t)\right|>0}
\]
and for $\varepsilon>0$ by $D(\varepsilon)$ we denote the set of
points where $x$ has jumps of size $\ge\varepsilon$, 
\[
D\rbr{\varepsilon}=\cbr{t\in(0,T]:\left|\Delta x(t)\right|\ge\varepsilon}.
\]
Using $D(\varepsilon)$ we define $J_{\varepsilon}(x):[0,T]\ra\R$
as
\[
J_{\varepsilon}(x)(t):=\sum_{s\in D(\varepsilon)\cap[0,t]}\Delta x(s).
\]

For $y:B\ra\R$ where $B\subseteq\R$, by the oscillation of $y$
on $A\subseteq B$ we will mean the quantity $\omega\rbr{y;A}$ defined
as 
\[
\omega\rbr{y;A}:=\sup_{s,t\in A}\left|y(t)-y(s)\right|
\]
and by the oscillation of $y$ on $A$ for increments not exceeding
$\varepsilon>0$ we will mean the quantity $\omega\rbr{y;\varepsilon;A}$
defined as 
\[
\omega\rbr{y;\varepsilon;A}:=\sup_{s,t\in A,|t-s|\le\varepsilon}\left|y(t)-y(s)\right|.
\]
Whenever $y$
is continuous, $A\neq \emptyset$ and $A$ is compact we naturally have $\lim_{\varepsilon\ra0+}\omega\rbr{y;\varepsilon;A}=0$.

Additionally, following \cite{HiraiEJP2023}, for a partition $\pi$ of the interval
$(0,T]$ we define 
\[
O(y;\pi):=\max_{(u,v]\in\pi}\omega\rbr{y;(u,v]}.
\]
(If $y$ is c\`adl\`ag, then also
\[
O(y;\pi)=\max_{(u,v]\in\pi}\omega\rbr{y;(u,v]}=\max_{(u,v]\in\pi}\omega\rbr{y;[u,v]}.
\]
)

We will consider sequences of partitions $\rbr{\pi^{n}}$ satisfying
the following two assumptions:

(A1)
\[
 \limsup_{\varepsilon\ra0+}\limsup_{n\ra+\ns}O(x-J_{\varepsilon}(x);\pi^{n})=0.
\]

(A2) For any $s\in(0,T]$, 
\[
\lim_{n\ra+\ns}x\rbr{\underline{\pi^{n}}(s)}=x(s-).
\]

Any sequence of partitions with vanishing mesh, that is satisfying \eqref{eq:mesh},  satisfies (A1) and (A2) as well, see \cite[Proposition 4.5]{HiraiEJP2023}. However, the opposite is not true. For example, the constant sequence consisting of two element partition $\pi=\cbr{(0,1/2],(1/2,1]}$  of $(0,1]$
for the indicator function ${\bf 1}_{[1/2,1]}:[0,1] \ra \cbr{0,1}$ of the interval $[1/2,1]$ satisfies (A1)-(A2) but not (\ref{eq:mesh}). 

\section{Main results}
In this section we present two main results of this note -- Theorem \ref{thm1} and Corollary \ref{cor1}. Their proofs are presented in the next section. 
\begin{theo} \label{thm1}
Let $x:[0,T]\ra\R$ be a c\`adl\`ag function. Let $\rbr{\pi^{n}}$
satisfies (A1) and (A2), and $[x]:[0,T]\ra[0,+\ns)$ be a non-decreasing
function such that the limit (\ref{eq:qvar1}) exists for all $t$
where $[x]$ is continuous, for $t=T$ and for all $t$ where $x$
has jumps, and at all these points
\[
\lim_{n\ra+\ns}\sum_{(u,v]\in\pi^{n}}\rbr{x\rbr{v\wedge t}-x\rbr{u\wedge t}}^{2}=[x]_{t}.
\]
Then for any twice continuously differentiable $f:\R\ra\R$, the finite limit $$I:=\lim_{n\ra+\ns}\sum_{(u,v]\in\pi^{n}}f''\rbr{x\rbr u}\rbr{x\rbr v-x\rbr u}^{2}$$
exists as well. As a result, the limit $\lim_{n\ra+\ns} \sum_{(u,v]\in\pi^{n}}f'\rbr{x\rbr u}\rbr{x\rbr v-x\rbr u}$ also
exists and is equal $f(x(T))-f(x(0))-\frac{1}{2}I-J$, where 
\begin{equation}
J=\sum_{t\in D}\cbr{f(x\rbr t)-f(x\rbr{t-})-f'\rbr{x\rbr{t-}}\Delta x(t)-\frac{1}{2}f''\rbr{x\rbr{t-}}\rbr{\Delta x\rbr t}^{2}}.\label{eq:jumps}
\end{equation}
\end{theo}
\begin{rema} Notice that we have not assumed that the jumps of $[x]$ occur
only at the points where $x$ has jumps, they also may occur at the
continuity points of $x$ as well. Next, we have not assumed anything
about the size of jumps. F\"ollmer assumed that 
\[
\Delta[x]_{t}:=[x]_{t}-[x]_{t-}=\rbr{\Delta x_{t}}^{2}=\rbr{x_{t}-x_{t-}}^{2}.
\]
In our case we always have $\Delta[x]_{t}\ge\rbr{\Delta x_{s}}^{2}$
(this follows from (A2)) but this inequality may be strict.
\end{rema}

\begin{rema} If the limit (\ref{eq:qvar1}) does not exist at some point
$t\in(0,T]$ where $x$ has a jump, then there may exist a twice continuously
differentiable function $f$ such that the limit $\lim_{n\ra+\ns}\sum_{(u,v]\in\pi^{n}}f''\rbr{x\rbr u}\rbr{x\rbr v-x\rbr u}^{2}$
does not exist either. 
\end{rema}
{\bf Example} Let us recall the function $z$ and
the Lebesgue partition $\rho^{n}(\alpha)$ of the interval $(0,1]$
($n\in\N$, $\alpha\in[0,1)$) for $z$. We define the zigzag function
$p:[0,2]\ra\R$ in the following way:
\[
p(t)=\begin{cases}
z(t) & \text{for }t\in[0,1),\\
2+z(2-t) & \text{for }t\in[1,2].
\end{cases}
\]
Let $\rho^{n}(\alpha)$ be equal the sequence of the intervals $\rbr{\left(t_{i}^{n}(\alpha),t_{i+1}^{n}(\alpha)\right]_{i=0}^{k_{n}(\alpha)-1}}$.
We define the partition $\sigma^{n}$, $n\in\N$, in the following
way: if $n$ is odd then $\sigma^{n}$ is the concatenation of the
sequences $\rho^{n}(0)$ and $\rbr{\left(2-t_{i+1}^{n}(1/2),2-t_{i}^{n}(1/2)\right]_{i=k_{n}(1/2)-1}^{0}}$
(this is a reflection of the sequence $\rho^{n}(1/2)$), while
if $n$ is even then $\sigma^{n}$ is the concatenation of the sequences
$\rho^{n}(1/2)$ and $\rbr{\left(2-t_{i+1}^{n}(0),2-t_{i}^{n}(0)\right]_{i=k_{n}(0)-1}^{0}}$
(this is a reflection of the sequence $\rho^{n}(0)$). From the
properties of the quadratic variation of $z$ along $\rbr{\rho^{n}(\alpha)}$
we get that 
\[
\lim_{n\ra+\ns}\sum_{(u,v]\in\sigma^{2n+1}}\rbr{p\rbr{v\wedge t}-p\rbr{u\wedge t}}^{2}=\begin{cases}
0 & \text{for }t\in[0,1),\\
L(0) +4 & \text{for }t=1,\\
L(0) + 4+L(1/2)& \text{for }t\in(1,2]
\end{cases}
\]
(this follows from the fact that for any $0\le  s < t < 1$ or $1< s <t \le 2$
\[
\lim_{n\ra+\ns} \cbr{ \sum_{(u,v]\in\sigma^{n}}\rbr{p\rbr{v\wedge t}-p\rbr{u\wedge t}}^{2} -\rbr{p\rbr{v\wedge s}-p\rbr{u\wedge s}}^{2}} = 0 
\]
since $p$ has finite total variation on $[s,t]$, moreover
\[
\lim_{n\ra+\ns}  \sum_{(u,v]\in\sigma^{2n+1}}\rbr{p\rbr{v\wedge 1}-p\rbr{u\wedge 1}}^{2} = \lim_{n\ra+\ns}  \sum_{(u,v]\in\rho^{2n+1}(0)} \rbr{z\rbr{v\wedge 1}-z\rbr{u\wedge 1}}^{2} + \rbr{\Delta p(1)}^2
\]
and 
\begin{align*}
\lim_{n\ra+\ns}  \sum_{(u,v]\in\sigma^{2n+1}}\rbr{p\rbr{v\wedge 2}-p\rbr{u\wedge 2}}^{2} = & \lim_{n\ra+\ns}  \sum_{(u,v]\in\rho^{2n+1}(0)} \rbr{z\rbr{v}-z\rbr{u}}^{2} + \rbr{\Delta p(1)}^2 \\
 & +  \lim_{n\ra+\ns}  \sum_{(u,v]\in\rho^{2n+1}(1/2)} \rbr{z\rbr{v}-z\rbr{u}}^{2} 
\end{align*}
)
while for even-numbered partitions we get 
\[
\lim_{n\ra+\ns}\sum_{(u,v]\in\sigma^{2n}}\rbr{p\rbr{v\wedge t}-p\rbr{u\wedge t}}^{2}=\begin{cases}
0 & \text{for }t\in[0,1),\\
L(1/2) +4 & \text{for }t=1,\\
L(1/2) + 4+L(0)& \text{for }t\in(1,2].
\end{cases}
\]
Thus, the quadratic variation of $p$ along $\rbr{\sigma^{n}}$ does
exist till any $t\in[0,2]$ excluding $t=1$. Now, if we have twice
continuously differentiable $f$ such that $f''(w)=1$ for $w\le1$
and $f''(w)=0$ for $w\ge2$ then 
\begin{align*}
\lim_{n\ra+\ns}\sum_{(u,v]\in\sigma^{n}}f''(p(u))\rbr{p\rbr v-p\rbr u}^{2} & =\lim_{n\ra+\ns}\sum_{(u,v]\in\sigma^{n},u\le1}f''(p(u))\rbr{p\rbr v-p\rbr u}^{2}\\
 & =\lim_{n\ra+\ns}\sum_{(u,v]\in\sigma^{n},u\le1}\rbr{p\rbr v-p\rbr u}^{2}\\
 & =\lim_{n\ra+\ns}\sum_{(u,v]\in\sigma^{n}}\rbr{p\rbr{v\wedge1}-p\rbr{u\wedge1}}^{2},
\end{align*}
but the last limit does not exist. 
\begin{coro} \label{cor1} Under the assumptions of Theorem \ref{thm1}, if $[x]$ is
right-continuous at all the points $t\in(0,T)$ where $x$ has jumps
then the limit $\lim_{n\ra+\ns}\sum_{(u,v]\in\pi^{n}}f''\rbr{x\rbr u}\rbr{x\rbr v-x\rbr u}^{2}$
may be expressed as the Lebesgue-Stieltjes integral:
\[
\lim_{n\ra+\ns}\sum_{(u,v]\in\pi^{n}}f''\rbr{x\rbr u}\rbr{x\rbr v-x\rbr u}^{2}=\int_{(0,T]}f''(x(t-))\dd{\widetilde{[x]}_{t}},
\]
where $\widetilde{[x]}$ denotes the right-continuous modification
of $[x]$, and we have the It{\^o}-F\"ollmer formula 
\[
f(x(T))-f(x(0))=\int_{(0,T]}f'(x(t-))\dd x(t)+\frac{1}{2}\int_{(0,T]}f''(x(t-))\dd{\widetilde{[x]}_{t}}+J,
\]
where by $\int_{(0,T]}f'(x(t-))\dd x(t)$ we denote the limit 
\[
\lim_{n\ra+\ns}\sum_{(u,v]\in\pi^{n}}f'\rbr{x\rbr u}\rbr{x\rbr v-x\rbr u}
\]
while $J$ is defined in (\ref{eq:jumps}). 
\end{coro}
\begin{rema}
If the limit (\ref{eq:qvar1}) is not right-continuous
at some point where $x$ has a jump then there may be no finite Borel
measure $\mu$ such that for any twice continuously differentiable
$f:\R\ra\R$, 
\[
\lim_{n\ra+\ns}\sum_{(u,v]\in\pi^{n}}f''\rbr{x\rbr u}\rbr{x\rbr v-x\rbr u}^{2}=\int_{(0,T]}f''(x(t-))\dd{\mu(t)}.
\]
\end{rema}
{\bf Example} Recall the function $z$. We define the zigzag function $q:[0,2]\ra\R$
in the following way:
\[
q(t)=\begin{cases}
-z(t)+t-1  & \text{for }t \in[0,1),\\
z(2-t)+t  & \text{for }t\in[1,2],
\end{cases}
\begin{cases}
\le 0 & \text{for }t\in[0,1),\\
\ge 1& \text{for }t\in[1,2].
\end{cases}
\]
Recall the partition $\rho^{n}(\alpha)$, $\alpha\in[0,1)$, and assume that it is equal the following sequence of intervals $\rbr{\left(t_{i}^{n}(\alpha),t_{i+1}^{n}(\alpha)\right]_{i=0}^{k_{n}(\alpha)-1}}$.
We define the partition $\tau^{n}$ of $(-1,1]$, $n\in\N$,  in the following way:
$\tau^{n}$ is the concatenation of the sequences $\rho^{n}(0)$ and
$\rbr{\left(2-t_{i+1}^{n}(0),2-t_{i}^{n}(0)\right]_{i=k_{n}(0)-1}^{0}}$.
Then the quadratic variation of $q$ along the sequence of partitions $\tau^n$, $n \in \N$, is not right continuous
at $1$, namely
\[
[q]_{1+}=[q]_{1}+L(0).
\]
Now let us consider functions $f_{m}$, $m\in\N$, such that $f_{m}''$
is equal $0$ on $(-\ns,0]$, it is equal $1$ on $[1,1+1/m]$, then 
$0$ on $[1+2/m,+\ns)$, $f_{m}''$ is affine on $[0,1]$ and on $[1+1/m,1+2/m]$.
For any $m\in\N$
\[
\lim_{n\ra+\ns}\sum_{(u,v]\in\pi^{n}}f''_{m}\rbr{q\rbr u}\rbr{q\rbr v-q\rbr u}^{2}\ge[q]_{1+}-[q]_{1}=L(0).
\]
If there was $\mu$ such that
\[
\lim_{n\ra+\ns}\sum_{(u,v]\in\pi^{n}}f''_{m}\rbr{q\rbr u}\rbr{q\rbr v-q\rbr u}^{2}=\int_{(0,2]}f''_{m}(q(t-))\dd{\mu(t)}
\]
then for any $m\in\N$,
\[
\int_{(0,2]}f''_{m}(q(t-))\dd{\mu(t)}\ge L(0)=\pi^{2}/3>0.
\]
On the other hand $f''_{m}(t)\ra t{\bf 1}_{[0,1]}(t)$,
thus the Lebesgue dominated convergence theorem yields
\[
\lim_{m\ra+\ns}\int_{(0,2]}f''_{m}(q(t-))\dd{\mu(t)}=\int_{(0,2]}q(t-){\bf 1}_{[0,1]}(q(t-))\dd{\mu(t)}=0.
\]
since for any $t\in(0,2]$, $q(t-)\in(-\ns,0]\cup(1,+\ns)$. 
\section{Proofs}
{\bf Proof of Theorem \ref{thm1}.} Let us fix $\varepsilon\in(0,1)$ and construct a
piecewise constant approximation of $x$, which approximates $x$
uniformly with accuracy $\varepsilon/3$. Our construction is the
standard one: let $t_{0}(\varepsilon):=0$ and having defined $t_{0}(\varepsilon)$,
..., $t_{k-1}(\varepsilon)$ such that $t_{k-1}(\varepsilon)<T$ for
some $k=1,2,\ldots$, we define
\begin{equation}
t_{k}(\varepsilon):=\inf\cbr{t>t_{k-1}(\varepsilon):\left|x(t)-x\rbr{t_{k-1}(\varepsilon)}\right|\ge\frac{1}{3}\varepsilon}\wedge T.\label{eq:tkdefinition}
\end{equation}
(in the above definition we apply the standard assumption that $\inf\emptyset=+\ns$).
Notice that by this definition and right continuity of $x$, if $t_{k}(\varepsilon)<T$
then
\[
\left|x\rbr{t_{k}(\varepsilon)}-x\rbr{t_{k-1}(\varepsilon)}\right|\ge\frac{1}{3}\varepsilon.
\]
We continue this procedure until $t_{k}(\varepsilon)=T$. We denote
the positive integer $k$ such that $t_{k}(\varepsilon)=T$ by $k(\varepsilon)$.
Let us also notice that 
\[
D(\varepsilon)\subseteq\cbr{t_{1}(\varepsilon),...,t_{k(\varepsilon)}(\varepsilon)}.
\]
Now, to deal with the situations when the limit (\ref{eq:qvar1})
does not exist, let $K(\varepsilon)=\cbr{0,\ldots,k(\varepsilon)}$
and let $L(\varepsilon)$ denote the subset of $K(\varepsilon)\setminus\cbr 0$
such that $x$ is continuous at $t_{l}$ whenever $l\in L(\varepsilon)$.
Recalling that, by assumption, the limit (\ref{eq:qvar1}) does exist
at all jump points of $x$, for any $l\in L(\varepsilon)$, $x$ is
continuous at $t_{l}(\varepsilon)$. Thus we may choose $t\in\rbr{t_{l-1}(\varepsilon),t_{l}(\varepsilon)}$
such that $x$ is continuous at $t$, 
\begin{equation}
\sup_{u\in\sbr{t,t_{l}(\varepsilon)}}\left|x\rbr{t_{l}(\varepsilon)}-x\rbr u\right|\le\frac{1}{3k(\varepsilon)}\varepsilon\label{eq:aaam}
\end{equation}
and $[x]$ is continuous at $t$ (both --
$x$ and $[x]$ have only countable many discontinuity points, $x$
because it is c\`adl\`ag while $[x]$ because it is non-decreasing). By
assumption, the limit (\ref{eq:qvar1}) does exist at
$t$. We denote such $t$ by $u_{l}(\varepsilon)$. For all $k\in K(\varepsilon)\setminus L(\varepsilon)$
we set $u_{k}(\varepsilon):=t_{k}(\varepsilon)$. Now we define
\begin{align*}
x^{\varepsilon}(t)= & \sum_{k=0}^{k(\varepsilon)-1}x\rbr{t_{k}(\varepsilon)}{\bf 1}_{\left[u_{k}(\varepsilon),u_{k+1}(\varepsilon)\right)}(t)+x\rbr{t_{k-1}(\varepsilon)}{\bf 1}_{[-\varepsilon/3,\varepsilon/3]}\rbr{x\rbr{t_{k-1}(\varepsilon)}-x\rbr T}{\bf 1}_{\cbr T}(t)\\
 & +x(T)\rbr{1-{\bf 1}_{[-\varepsilon/3,\varepsilon/3]}\rbr{x\rbr{t_{k-1}(\varepsilon)}-x\rbr T}}{\bf 1}_{\cbr T}(t).
\end{align*}
(${\bf 1}_{A}$ is the indicator function of the set $A$.) Let us
make few observations about $x^{\varepsilon}$. $x^{\varepsilon}$
is c\`adl\`ag, it approximates $x$ uniformly with accuracy $\varepsilon/3$
and any jump of $x$ of size $\ge\varepsilon$ is also a jump of $x^{\varepsilon}$.

Denote $M=\sup_{u\in[0,T]}|x(u)|$. Now we estimate the difference
\begin{align}
 & \left|\sum_{(u,v]\in\pi^{n}}f''\rbr{x\rbr u}\rbr{x\rbr v-x\rbr u}^{2}-\sum_{(u,v]\in\pi^{n}}f''\rbr{x^{\varepsilon}\rbr u}\rbr{x\rbr v-x\rbr u}^{2}\right|\nonumber \\
 & \le\sum_{(u,v]\in\pi^{n}}\left|f''\rbr{x\rbr u}-f''\rbr{x^{\varepsilon}\rbr u}\right|\rbr{x\rbr v-x\rbr u}^{2}\nonumber \\
 & \le D(\varepsilon)\sum_{(u,v]\in\pi^{n}}\rbr{x\rbr v-x\rbr u}^{2},\label{eq:fbisdiffestim}
\end{align}
where 
\[
D(\varepsilon):=\omega\rbr{f'';\frac{1}{3}\varepsilon;\sbr{-M-\frac{1}{3},M+\frac{1}{3}}}.
\]
Recalling that $K(\varepsilon)=\cbr{0,\ldots,k(\varepsilon)}$ we
write 
\begin{align}
 \sum_{(u,v]\in\pi^{n}}f''\rbr{x^{\varepsilon}\rbr u}\rbr{x\rbr v-x\rbr u}^{2} = & \sum_{k=0}^{k(\varepsilon)-1}f''\rbr{x\rbr{t_{k}(\varepsilon)}}\sum_{(u,v]\in\pi^{n}:u_{k}(\varepsilon)\le u<v<u_{k+1}(\varepsilon)}\rbr{x\rbr v-x\rbr u}^{2}\nonumber \\
 & +\sum_{k=1}^{k(\varepsilon)}f''\rbr{x\rbr{t_{k-1}(\varepsilon)}}\sum_{(u,v]\in\pi^{n}:u<u_{k}(\varepsilon)\le v}\rbr{x\rbr v-x\rbr u}^{2}.\label{eq:expansion}
\end{align}
Now we are going to prove that for $k=1,2,\ldots,k(\varepsilon)-1$
\begin{equation}
\liminf_{n\ra+\ns}\sum_{(u,v]\in\pi^{n}:u_{k}(\varepsilon)\le u<v<u_{k+1}(\varepsilon)}\rbr{x\rbr v-x\rbr u}^{2}\ge[x]_{u_{k+1}(\varepsilon)}-\rbr{\Delta x\rbr{u_{k+1}(\varepsilon)}}^{2}-[x]_{u_{k}(\varepsilon)}-\frac{4M}{k(\varepsilon)}\varepsilon,\label{eq:cccm}
\end{equation}
\begin{equation}
\limsup_{n\ra+\ns}\sum_{(u,v]\in\pi^{n}:u_{k}(\varepsilon)\le u<v<u_{k+1}(\varepsilon)}\rbr{x\rbr v-x\rbr u}^{2}\le[x]_{u_{k+1}(\varepsilon)}-\rbr{\Delta x\rbr{u_{k+1}(\varepsilon)}}^{2}-[x]_{u_{k}(\varepsilon)}+\frac{4M}{k(\varepsilon)}\varepsilon\label{eq:cccm1}
\end{equation}
and for $k=1,2,\ldots,k(\varepsilon)$
\begin{equation}
\liminf_{n\ra+\ns}\sum_{(u,v]\in\pi^{n}:u<u_{k}(\varepsilon)\le v}\rbr{x\rbr v-x\rbr u}^{2}\ge\rbr{\Delta x\rbr{u_{k}(\varepsilon)}}^{2}-\frac{4M}{3k(\varepsilon)}\varepsilon,\label{eq:bbbm}
\end{equation}
\begin{equation}
\limsup_{n\ra+\ns}\sum_{(u,v]\in\pi^{n}:u<u_{k}(\varepsilon)\le v}\rbr{x\rbr v-x\rbr u}^{2}\ge\rbr{\Delta x\rbr{u_{k}(\varepsilon)}}^{2}+\frac{4M}{3k(\varepsilon)}\varepsilon,\label{eq:bbbm1}
\end{equation}
where $M:=\sup_{u\in[0,T]}|x(u)|$. Notice that by the choice of $u_{k}(\varepsilon)$,
$k\in K(\varepsilon)$, the limit (\ref{eq:qvar1}) exists for all
$t\in\cbr{u_{0}(\varepsilon),u_{1}(\varepsilon),\ldots,u_{k(\varepsilon)}(\varepsilon)}$. 

To prove (\ref{eq:cccm})-(\ref{eq:bbbm1}) first we observe that
\begin{align*}
\sum_{(u,v]\in\pi^{n}:u_{k}(\varepsilon)\le u<v<u_{k+1}(\varepsilon)}\rbr{x\rbr v-x\rbr u}^{2} &  =\sum_{(u,v]\in\pi^{n}:0\le u<v<u_{k+1}(\varepsilon)}\rbr{x\rbr v-x\rbr u}^{2} \\
& \quad -\sum_{(u,v]\in\pi^{n}:0\le u<v<u_{k}(\varepsilon)}\rbr{x\rbr v-x\rbr u}^{2} \\
 & \quad-\sum_{(u,v]\in\pi^{n}:u<u_{k}(\varepsilon)\le v}\rbr{x\rbr v-x\rbr u}^{2}
\end{align*}
hence to prove (\ref{eq:cccm}) it suffices to prove that for $k=1,2,\ldots,k(\varepsilon)$
\begin{equation}
\limsup_{n\ra+\ns}\rbr{\sum_{(u,v]\in\pi^{n}:0\le u<v<u_{k}(\varepsilon)}\rbr{x\rbr v-x\rbr u}^{2}+\sum_{(u,v]\in\pi^{n}:u<u_{k}(\varepsilon)\le v}\rbr{x\rbr v-x\rbr u}^{2}}\le[x]_{u_{k}(\varepsilon)}+\frac{4M}{3k(\varepsilon)}\varepsilon,\label{eq:dddm}
\end{equation}
\begin{equation}
\liminf_{n\ra+\ns}\rbr{\sum_{(u,v]\in\pi^{n}:0\le u<v<u_{k}(\varepsilon)}\rbr{x\rbr v-x\rbr u}^{2}+\sum_{(u,v]\in\pi^{n}:u<u_{k}(\varepsilon)\le v}\rbr{x\rbr v-x\rbr u}^{2}}\ge[x]_{u_{k}(\varepsilon)}-\frac{4M}{3k(\varepsilon)}\varepsilon\label{eq:dddm1}
\end{equation}
and (\ref{eq:bbbm})-(\ref{eq:bbbm1}).

To prove (\ref{eq:dddm})-(\ref{eq:dddm1}) let us fix $k=1,2,\ldots,k(\varepsilon)$
and observe 
\begin{align}
 & \left|\sum_{(u,v]\in\pi^{n}:0\le u<v<u_{k}(\varepsilon)}\rbr{x\rbr v-x\rbr u}^{2}+\sum_{(u,v]\in\pi^{n}:u<u_{k}(\varepsilon)\le v}\rbr{x\rbr v-x\rbr u}^{2}-\sum_{(u,v]\in\pi^{n}}\rbr{x\rbr{v\wedge u_{k}(\varepsilon)}-x\rbr{u\wedge u_{k}(\varepsilon)}}^{2}\right|\nonumber \\
 & =\left|\sum_{(u,v]\in\pi^{n}:u<u_{k}(\varepsilon)\le v}\rbr{x\rbr v-x\rbr u}^{2}-\sum_{(u,v]\in\pi^{n}:u<u_{k}(\varepsilon)\le v}\rbr{x\rbr{u_{k}(\varepsilon)}-x\rbr u}^{2}\right|\nonumber \\
 & =\left|\rbr{x\rbr{\overline{\pi^{n}}\rbr{u_{k}(\varepsilon)}}-x\rbr{\underline{\pi^{n}}\rbr{u_{k}(\varepsilon)}}}^{2}-\rbr{x\rbr{u_{k}(\varepsilon)}-x\rbr{\underline{\pi^{n}}\rbr{u_{k}(\varepsilon)}}}^{2}\right|\nonumber \\
 & =\left|x\rbr{\overline{\pi^{n}}\rbr{u_{k}(\varepsilon)}}+x\rbr{u_{k}(\varepsilon)}-2x\rbr{\underline{\pi^{n}}\rbr{u_{k}(\varepsilon)}}\right|\left|x\rbr{\overline{\pi^{n}}\rbr{u_{k}(\varepsilon)}}-x\rbr{u_{k}(\varepsilon)}\right|\nonumber \\
 & \le4M\left|x\rbr{\overline{\pi^{n}}\rbr{u_{k}(\varepsilon)}}-x\rbr{u_{k}(\varepsilon)}\right|.\label{eq:difference}
\end{align}

If $k=k(\varepsilon)$ then $u_{k}(\varepsilon)=T$ and then 
\begin{equation}
x\rbr{\overline{\pi^{n}}\rbr{u_{k}(\varepsilon)}}-x\rbr{u_{k}(\varepsilon)}=x(T)-x(T)=0.\label{eq:brr1}
\end{equation}

If $k<k(\varepsilon)$ and $u_{k}(\varepsilon)$ is a point of continuity
of $x$ then by (\ref{eq:aaam}) $x$ is continuous at $t_{k}(\varepsilon)$,
$t_{k}(\varepsilon)\in\left(u_{k}(\varepsilon),T\right]$ and
\begin{equation}
\sup_{u\in\sbr{u_{k}(\varepsilon),t_{k}(\varepsilon)}}\left|x\rbr u-x\rbr{t_{k}(\varepsilon)}\right|\le\frac{1}{3k(\varepsilon)}\varepsilon.\label{eq:estimuktk}
\end{equation}
Let us notice that $x\rbr{t_{k}(\varepsilon)}\neq x\rbr{u_{k}(\varepsilon)}$,
otherwise the infimum of the set on the right side of (\ref{eq:tkdefinition})
would occurr already at $u_{k}(\varepsilon)$ not at $t_{k}(\varepsilon)$.
Since $x\rbr{t_{k}(\varepsilon)}\neq x\rbr{u_{k}(\varepsilon)}$ and
$x$ is continuous at $u_{k}(\varepsilon)$ and $t_{k}(\varepsilon)$,
by (A2), for sufficiently large $n\in\N$ we have 
\[
\underline{\pi^{n}}\rbr{u_{k}(\varepsilon)}<u_{k}(\varepsilon)\le\overline{\pi^{n}}\rbr{u_{k}(\varepsilon)}\le\underline{\pi^{n}}\rbr{t_{k}(\varepsilon)}<t_{k}(\varepsilon)
\]
 and by (\ref{eq:estimuktk}) 
\begin{equation}
\limsup_{n\ra+\ns}\left|x\rbr{\overline{\pi^{n}}\rbr{u_{k}(\varepsilon)}}-x\rbr{u_{k}(\varepsilon)}\right|\le\frac{1}{3k(\varepsilon)}\varepsilon.\label{eq:brr2}
\end{equation}

Let now $k<k(\varepsilon)$ and $u_{k}(\varepsilon)$ be a point where
$x$ has a jump. By right continuity of $x$, there exists some $v\in\rbr{u_{k}(\varepsilon),T}$
such that $x$ is continuous at $v$, $x\rbr{u_{k}(\varepsilon)-}\neq x\rbr v$
and
\begin{equation}
\sup_{u\in\sbr{u_{k}(\varepsilon),v}}\left|x\rbr u-x\rbr{u_{k}(\varepsilon)}\right|\le\frac{1}{3k(\varepsilon)}\varepsilon\label{eq:estimukv}
\end{equation}
Now, again by (A2), for sufficiently large $n\in\N$ we have 
\[
\underline{\pi^{n}}\rbr{u_{k}(\varepsilon)}<u_{k}(\varepsilon)\le\overline{\pi^{n}}\rbr{u_{k}(\varepsilon)}\le\underline{\pi^{n}}\rbr v<v
\]
and by (\ref{eq:estimukv}) 
\begin{equation}
\limsup_{n\ra+\ns}\left|x\rbr{\overline{\pi^{n}}\rbr{u_{k}(\varepsilon)}}-x\rbr{u_{k}(\varepsilon)}\right|\le\frac{1}{3k(\varepsilon)}\varepsilon.\label{eq:brr3}
\end{equation}
Now, using (\ref{eq:difference}) and (\ref{eq:brr1}), (\ref{eq:brr2})
or (\ref{eq:brr3}) we obtain (\ref{eq:dddm})-(\ref{eq:dddm1}). 

To prove (\ref{eq:bbbm})-(\ref{eq:bbbm1}) we again use the estimate
of 
\[
\left|\sum_{(u,v]\in\pi^{n}:u<u_{k}(\varepsilon)\le v}\rbr{x\rbr v-x\rbr u}^{2}-\sum_{(u,v]\in\pi^{n}:u<u_{k}(\varepsilon)\le v}\rbr{x\rbr{u_{k}(\varepsilon)}-x\rbr u}^{2}\right|
\]
obtained in (\ref{eq:difference}), then (\ref{eq:brr1}), (\ref{eq:brr2})
or (\ref{eq:brr3}) and the fact that by (A2)
\begin{align*}
 & \lim_{n\ra+\ns}\sum_{(u,v]\in\pi^{n}:u<u_{k}(\varepsilon)\le v}\rbr{x\rbr{u_{k}(\varepsilon)}-x\rbr u}^{2}\\
 & =\lim_{n\ra+\ns}\rbr{x\rbr{u_{k}(\varepsilon)}-x\rbr{\underline{\pi^{n}}\rbr{u_{k}(\varepsilon)}}}^{2}=\rbr{\Delta x\rbr{u_{k}(\varepsilon)}}^{2}.
\end{align*}

Let us denote
\begin{align}
I(\varepsilon): & =\sum_{k=0}^{k(\varepsilon)-1}f''\rbr{x\rbr{t_{k}(\varepsilon)}}\cbr{[x]_{u_{k+1}(\varepsilon)}-\rbr{\Delta x\rbr{u_{k+1}(\varepsilon)}}^{2}-[x]_{u_{k}(\varepsilon)}}+\sum_{k=1}^{k(\varepsilon)}f''\rbr{x\rbr{t_{k-1}(\varepsilon)}}\rbr{\Delta x\rbr{u_{k}(\varepsilon)}}^{2}\nonumber \\
 & =\sum_{k=0}^{k(\varepsilon)-1}f''\rbr{x\rbr{t_{k}(\varepsilon)}}\cbr{[x]_{u_{k+1}(\varepsilon)}-[x]_{u_{k}(\varepsilon)}}.\label{eq:Ierepresentation}
\end{align}
Now, using estimates (\ref{eq:cccm})-(\ref{eq:bbbm1}) and (\ref{eq:expansion})
we get:
\begin{align*}
\limsup_{n\ra+\ns}\sum_{(u,v]\in\pi^{n}}f''\rbr{x^{\varepsilon}\rbr u}\rbr{x\rbr v-x\rbr u}^{2}-I(\varepsilon) & \le\sum_{k=1}^{k(\varepsilon)-1}\left|f''\rbr{x\rbr{t_{k}(\varepsilon)}}\right|\frac{4M}{k(\varepsilon)}\varepsilon+\sum_{k=1}^{k(\varepsilon)}f''\rbr{x\rbr{t_{k-1}(\varepsilon)}}\frac{4M}{3k(\varepsilon)}\varepsilon\\
 & \le k(\varepsilon)\omega\rbr{f'';[-M,M]}\rbr{\frac{4M}{k(\varepsilon)}+\frac{4M}{3k(\varepsilon)}}\varepsilon\\
 & =\frac{16}{3}\omega\rbr{f'';[-M,M]}M\varepsilon
\end{align*}
and similarly, 
\[
\liminf_{n\ra+\ns}\sum_{(u,v]\in\pi^{n}}f''\rbr{x^{\varepsilon}\rbr u}\rbr{x\rbr v-x\rbr u}^{2}-I(\varepsilon)\ge-\frac{16}{3}\omega\rbr{f'';[-M,M]}M\varepsilon.
\]
These together with the bound (\ref{eq:fbisdiffestim}) yield 
\[
\limsup_{n\ra+\ns}\left|\sum_{(u,v]\in\pi^{n}}f''\rbr{x\rbr u}\rbr{x\rbr v-x\rbr u}^{2}-I(\varepsilon)\right|\le[x]_{T}D(\varepsilon)+\frac{16}{3}\omega\rbr{f'';[-M,M]}M\varepsilon.
\]
Since $\sum_{(u,v]\in\pi^{n}}f''\rbr{x\rbr u}\rbr{x\rbr v-x\rbr u}^{2}$
does not depend on $\varepsilon$ and $D(\varepsilon)\ra0$ as $\varepsilon\ra0+$
we get that $I\rbr{\varepsilon_{1}}$ and $I\rbr{\varepsilon_{2}}$,
$\varepsilon_{1},\varepsilon_{2}\in(0,1)$ are as close as we please,
provided $\varepsilon_{1}$ and $\varepsilon_{2}$ are sufficiently
close to $0$, hence $\lim_{\varepsilon\ra0+}I(\varepsilon)$ exists
and thus the limit
\[
\lim_{n\ra+\ns}\sum_{(u,v]\in\pi^{n}}f''\rbr{x\rbr u}\rbr{x\rbr v-x\rbr u}^{2}
\]
exists and they are equal.

Now we are going to prove that the limit 
\[
\sum_{(u,v]\in\pi^{n}}f'\rbr{x\rbr u}\rbr{x\rbr v-x\rbr u}
\]
exists as well. We use the standard approach and write 
\[
f(x(T))-f(x(0))-\sum_{(u,v]\in\pi^{n}}f'\rbr{x\rbr u}\rbr{x\rbr v-x\rbr u}=\sum_{(u,v]\in\pi^{n}}\cbr{f(x\rbr v)-f(x\rbr u)-f'\rbr{x\rbr u}\rbr{x\rbr v-x\rbr u}}.
\]
Next, we use Taylor's expansion in the form 
\begin{align}
f(z)-f(y)-f'(y)(z-u) & =\int_{y}^{z}f''(w)(z-w)\dd w\nonumber \\
 & =\frac{1}{2}f''(y)(z-y)^{2}+\int_{y}^{z}\cbr{f''(w)-f''(y)}(z-w)\dd w\nonumber \\
 & =\frac{1}{2}f''(y)(z-y)^{2}+R(y,z)(z-y)^{2},\label{eq:Taylor}
\end{align}
where 
\begin{equation}
\left|R(y,z)\right|\le\frac{1}{2}\sup_{w\in[y\wedge z,y\vee z]}\left|f''\rbr w-f''(y)\right|.\label{eq:Restim}
\end{equation}
Using similar expansion as (\ref{eq:expansion}) and (\ref{eq:Taylor})
we get 
\begin{align*}
 & f(x(T))-f(x(0))-\sum_{(u,v]\in\pi^{n}}f'\rbr{x\rbr u}\rbr{x\rbr v-x\rbr u}-\frac{1}{2}\sum_{(u,v]\in\pi^{n}}f''\rbr{x\rbr u}\rbr{x\rbr v-x\rbr u}^{2}\\
 & =\sum_{(u,v]\in\pi^{n}}\cbr{f(x\rbr v)-f(x\rbr u)-f'\rbr{x\rbr u}\rbr{x\rbr v-x\rbr u}-\frac{1}{2}f''\rbr{x\rbr u}\rbr{x\rbr v-x\rbr u}^{2}}\\
 & =\sum_{k=0}^{k(\varepsilon)-1}\sum_{(u,v]\in\pi^{n}:u_{k}(\varepsilon)\le u<v<u_{k+1}(\varepsilon)}R\rbr{x\rbr u,x\rbr v}\rbr{x\rbr v-x\rbr u}^{2}\\
 & \quad+\sum_{k=1}^{k(\varepsilon)}\sum_{(u,v]\in\pi^{n}:u_{k}(\varepsilon)\in D(\varepsilon),u<u_{k}(\varepsilon)\le v}\cbr{f(x\rbr v)-f(x\rbr u)-f'\rbr{x\rbr u}\rbr{x\rbr v-x\rbr u}-\frac{1}{2}f''\rbr{x\rbr u}\rbr{x\rbr v-x\rbr u}^{2}}\\
 & \quad+\sum_{k=1}^{k(\varepsilon)}\sum_{(u,v]\in\pi^{n}:u_{k}(\varepsilon)\notin D(\varepsilon),u<u_{k}(\varepsilon)\le v}R\rbr{x\rbr u,x\rbr v}\rbr{x\rbr v-x\rbr u}^{2}.
\end{align*}
Let us notice that for $(u,v]\in\pi^{n}$ such that $u_{k}(\varepsilon)\le u<v<u_{k+1}(\varepsilon)$
for some $k=0,1,\ldots,k(\varepsilon)-1$, $(u,v]\cap D(\varepsilon)=\emptyset$
hence $O\rbr{x;\pi^{n}}=O\rbr{x-J_{\varepsilon}(x);\pi^{n}}$ and,
using also (\ref{eq:Restim}) we get
\begin{align*}
 & \limsup_{n\ra+\ns}\sum_{k=0}^{k(\varepsilon)-1}\sum_{(u,v]\in\pi^{n}:u_{k}(\varepsilon)\le u<v<u_{k+1}(\varepsilon)}\left|R\rbr{x\rbr u,x\rbr v}\right|\rbr{x\rbr v-x\rbr u}^{2}\\
 & \le\frac{1}{2}\limsup_{n\ra+\ns}\omega\rbr{f'';O\rbr{x-J_{\varepsilon}(x);\pi^{n}};[-M,M]}[x]_{T}.
\end{align*}
Next, using similar reasoning as in the proof of
(\ref{eq:brr2}) and (\ref{eq:brr3}), and (A2), for $u_{k}(\varepsilon)\notin D(\varepsilon)$
we have 
\[
\limsup_{n\ra+\ns}\sum_{(u,v]\in\pi^{n}:u_{k}(\varepsilon)\notin D(\varepsilon),u<u_{k}(\varepsilon)\le v}\omega\rbr{x;(u,v]}\le\varepsilon
\]
hence, by (\ref{eq:Restim})
\[
\limsup_{n\ra+\ns}\sum_{(u,v]\in\pi^{n}:u_{k}(\varepsilon)\notin D(\varepsilon),u<u_{k}(\varepsilon)\le v}\left|R\rbr{x\rbr u,x\rbr v}\right|\le\omega\rbr{f'';\varepsilon;[-M,M]}
\]
and
\begin{align*}
 & \limsup_{n\ra+\ns}\sum_{k=1}^{k(\varepsilon)}\sum_{(u,v]\in\pi^{n}:u_{k}(\varepsilon)\notin D(\varepsilon),u<u_{k}(\varepsilon)\le v}\left|R\rbr{x\rbr u,x\rbr v}\right|\rbr{x\rbr v-x\rbr u}^{2}\\
 & \le\frac{1}{2}\limsup_{n\ra+\ns}\omega\rbr{f'';\varepsilon;[-M,M]}[x]_{T}.
\end{align*}
Using the same argument which was used to prove (\ref{eq:brr3}) and
the fact that there are only finitely many jumps of $x$ of size $\ge\varepsilon$,
we also have 
\begin{align*}
 & \lim_{n\ra+\ns}\sum_{k=1}^{k(\varepsilon)}\sum_{(u,v]\in\pi^{n}:u_{k}(\varepsilon)\in D(\varepsilon),u<u_{k}(\varepsilon)\le v}\cbr{f(x\rbr v)-f(x\rbr u)-f'\rbr{x\rbr u}\rbr{x\rbr v-x\rbr u}-\frac{1}{2}f''\rbr{x\rbr u}\rbr{x\rbr v-x\rbr u}^{2}}\\
 & =\sum_{t\in D(\varepsilon)}\cbr{f(x\rbr t)-f(x\rbr{t-})-f'\rbr{x\rbr{t-}}\rbr{x\rbr t-x\rbr{t-}}-\frac{1}{2}f''\rbr{x\rbr{t-}}\rbr{x\rbr t-x\rbr{t-}}^{2}}\\
 & =\sum_{t\in D(\varepsilon)}\cbr{\Delta f(x\rbr t)-f'\rbr{x\rbr{t-}}\Delta x(t)-\frac{1}{2}f''\rbr{x\rbr{t-}}\rbr{\Delta x\rbr t}^{2}}.
\end{align*}
As a result, denoting 
\[
I=\lim_{n\ra+\ns}\sum_{(u,v]\in\pi^{n}}f''\rbr{x\rbr u}\rbr{x\rbr v-x\rbr u}^{2}
\]
 and 
\[
J=\sum_{t\in D}\cbr{\Delta f(x\rbr t)-f'\rbr{x\rbr{t-}}\Delta x(t)-\frac{1}{2}f''\rbr{x\rbr{t-}}\rbr{\Delta x\rbr t}^{2}},
\]
\[
J(\varepsilon)=\sum_{t\in D(\varepsilon)}\cbr{\Delta f(x\rbr t)-f'\rbr{x\rbr{t-}}\Delta x(t)-\frac{1}{2}f''\rbr{x\rbr{t-}}\rbr{\Delta x\rbr t}^{2}}
\]
we get 
\begin{align*}
 & \limsup_{n\ra+\ns}\left|\sum_{(u,v]\in\pi^{n}}f'\rbr{x\rbr u}\rbr{x\rbr v-x\rbr u}-\cbr{f(x(T))-f(x(0))-\frac{1}{2}I-J}\right|\\
 & \le\limsup_{n\ra+\ns}\omega\rbr{f'';O\rbr{x-J_{\varepsilon}(x);\pi^{n}}+\varepsilon;[-M,M]}[x]_{T}+J-J(\varepsilon).
\end{align*}
Since 
\[
\limsup_{\varepsilon\ra0}\limsup_{n\ra+\ns}\omega\rbr{f'';O\rbr{x-J_{\varepsilon}(x);\pi^{n}}+\varepsilon;[-M,M]}=0,
\]
$\limsup_{\varepsilon\ra0+}\rbr{J-J(\varepsilon)}=0$ and $\varepsilon$
may be chosen as small as we please, we get that the limit $$\lim_{n\ra+\ns}\sum_{(u,v]\in\pi^{n}}f'\rbr{x\rbr u}\rbr{x\rbr v-x\rbr u}$$ exists and is equal $f(x(T))-f(x(0))-\frac{1}{2}I-J$.
\hfill $\blacksquare$

\noindent {\bf Proof of Corollary \ref{cor1}.} We will use the same notation as in the proof of
Theorem \ref{thm1}.

Let $\widetilde{[x]}$ be a right-continuous modification of $t\mapsto[x]_{t}$.
Denote by $\tilde{\mu}$ the Lebesgue-Stieltjes measure induced by
$\widetilde{[x]}$, that means $\mu(a,b]:=\widetilde{[x]}_{b}-\widetilde{[x]}_{a}$,
and for a Borel-measurable $g:(0,T]\ra\R$ by $\int_{(0,T]}g(t)\dd{\widetilde{[x]}_{t}}$
we will denote the integral $\int_{(0,T]}g\dd{\tilde{\mu}}.$ Recall
the construction of $x^{\varepsilon}$. We notice that $\widetilde{[x]}_{t}=[x]_{t}$
for all $t\in\cbr{t_{0}(\varepsilon),t_{1}(\varepsilon),...,t_{k(\varepsilon)}(\varepsilon)}$.
Using the representation (\ref{eq:Ierepresentation}) and the fact
that by the construction $\left|x\rbr{u_{k}(\varepsilon)}-x\rbr{t_{k}(\varepsilon)}\right|\le\frac{1}{3}\varepsilon$,
$k\in K(\varepsilon)$, (recall (\ref{eq:aaam})) we see that we have
the estimate 
\begin{align*}
\left|I(\varepsilon)-\int_{(0,T]}f''(x^{\varepsilon}(t-))\dd{\widetilde{[x]}_{t}}\right| & \le\sum_{k=0}^{k(\varepsilon)-1}\left|f''\rbr{x\rbr{t_{k}(\varepsilon)}}-f''\rbr{x\rbr{u_{k}(\varepsilon)}}\right|\cbr{[x]_{u_{k+1}(\varepsilon)}-[x]_{u_{k}(\varepsilon)}}\\
 & \le\frac{1}{3}\omega\rbr{f'',\frac{1}{3}\varepsilon;\sbr{-M-\frac{1}{3},M+\frac{1}{3}}}[x]_{T}.
\end{align*}
By the construction, $x^{\varepsilon}\ra x$ uniformly as $\varepsilon\ra0+$.
This and the Lebesgue dominated convergence theorem yield
\[
\lim_{\varepsilon\ra0+}\int_{(0,T]}f''(x^{\varepsilon}(t-))\dd{\widetilde{[x]}_{t}}=\int_{(0,T]}f''(x(t-))\dd{\widetilde{[x]}_{t}}.
\]
On the other hand, $\lim_{\varepsilon\ra0+}I(\varepsilon)=I$, thus
$I=\int_{(0,T]}f''(x(t-))\dd{\widetilde{[x]}_{t}}$. 
\hfill $\blacksquare$

\section*{Appendix}
{\bf Proof of \eqref{zigzag}.}
We turn our attention to $k\rbr{1/\sqrt{n},\alpha/\sqrt{n}}$ -- the number of intervals in the Lebesgue partition
$\rho^{n}(\alpha) = {\pi\rbr{1/\sqrt{n},\alpha/\sqrt{n}}}$ for $z$, $n \in \N$, $\alpha \in [0, 1)$.

When $\alpha=0$ (the shift is $0$) then 
\begin{align*}
k\rbr{1/\sqrt{n},0} = 2\sum^\infty_{m=1}\left\lfloor\frac{\frac{1}{\sqrt{m}}}{\frac{1}{\sqrt{n}}}\right\rfloor.
\end{align*}

In the general case $\alpha\in[0,1)$ we have an almost similar formula as in the case of no shifts,
\begin{align*}
k\rbr{1/\sqrt{n},\alpha/\sqrt{n}}= 2\sum^\infty_{m=1}\left\lfloor\frac{\frac{1}{\sqrt{m}}-\frac{\alpha}{\sqrt{n}}}{\frac{1}{\sqrt{n}}}\right\rfloor
=2\sum^n_{m=1}\left\lfloor\sqrt{\frac{n}{m}}-\alpha\right\rfloor
\end{align*}
and
\[
\sum_{(u,v]\in\rho^{n}(\alpha)}  \rbr{z\rbr{v}-z\rbr{u}}^{2}  = \frac{1}{n} k\rbr{1/\sqrt{n},\alpha/\sqrt{n}}  =  \frac{2}{n} \sum^n_{m=1}\left\lfloor\sqrt{\frac{n}{m}}-\alpha\right\rfloor.
\]
Let
\begin{align*}
L_l:=\#\left\lbrace m\in \N: \left\lfloor\sqrt{\frac{n}{m}}-\alpha\right\rfloor=l\right\rbrace.
\end{align*}
We obtain that
\begin{align*}
k\rbr{1/\sqrt{n},\alpha/\sqrt{n}} = 2 \sum^n_{m=1}\left\lfloor\sqrt{\frac{n}{m}}-\alpha\right\rfloor= 2 \sum^{\lfloor\sqrt{n}\rfloor}_{l=1}l\cdot L_l.
\end{align*}
But $\left\lfloor\sqrt{\frac{n}{m}}-\alpha\right\rfloor=l$ if an only if $l\leq\sqrt{\frac{n}{m}}-\alpha<l+1$, and so
\begin{align*}
l+\alpha\leq\sqrt{\frac{n}{m}}<l+1+\alpha,
\end{align*}
and consequently,
\begin{align*}
\frac{n}{(l+1+\alpha)^2}<m\leq\frac{n}{(l+\alpha)^2}.
\end{align*}
Recall that for any real number $a$ we have $a=\lfloor a\rfloor+\{a\}$, where $\lfloor a\rfloor$ is floor of $a$ and $\{a\}$ is fractional part of $a$. Now, the bound on $L_l$ will be the number of integers between the real numbers $\frac{n}{(l+1+\alpha)^2}$ and $n\leq\frac{n}{(l+\alpha)^2}$, so $L_l$ is given by
\begin{align*}
L_l&=\left\lfloor\frac{n}{(l+\alpha)^2}\right\rfloor-\left\lfloor\frac{n}{(l+1+\alpha)^2}\right\rfloor\\
&=\frac{n}{(l+\alpha)^2}-\left\lbrace\frac{n}{(l+\alpha)^2}\right\rbrace-\frac{n}{(l+1+\alpha)^2}+\left\lbrace\frac{n}{(l+1+\alpha)^2}\right\rbrace.
\end{align*}
Now, we have
\begin{align*}
\begin{split}
\sum_{(u,v]\in\rho^{n}(\alpha)} & \rbr{z\rbr{v}-z\rbr{u}}^{2} = \frac{1}{n} k\rbr{1/\sqrt{n},\alpha/\sqrt{n}}  = \frac{2}{n}\sum^{\lfloor\sqrt{n}\rfloor}_{l=1}l\cdot L_l\\
=&\frac{2}{n}\sum^{\lfloor\sqrt{n}\rfloor}_{l=1}\left(\frac{n}{(l+\alpha)^2}-\frac{n}{(l+1+\alpha)^2}\right)l+\frac{2}{n}\sum^{\lfloor\sqrt{n}\rfloor}_{l=1}\left\lbrace\frac{n}{(l+1+\alpha)^2}\right\rbrace l-\frac{2}{n}\sum^{\lfloor\sqrt{n}\rfloor}_{l=1}\left\lbrace\frac{n}{(l+\alpha)^2}\right\rbrace l\\
=&\frac{2}{n}\sum^{\lfloor\sqrt{n}\rfloor}_{l=1}\left(\frac{n}{(l+\alpha)^2}-\frac{n}{(l+1+\alpha)^2}\right)l+\frac{2}{n}\sum^{\lfloor\sqrt{n}\rfloor}_{l=1}\left\lbrace\frac{n}{(l+1+\alpha)^2}\right\rbrace (l+1)\\
&-\frac{2}{n}\sum^{\lfloor\sqrt{n}\rfloor}_{l=1}\left\lbrace\frac{n}{(l+1+\alpha)^2}\right\rbrace-\frac{2}{n}\sum^{\lfloor\sqrt{n}\rfloor}_{l=1}\left\lbrace\frac{n}{(l+\alpha)^2}\right\rbrace l\\
=&\frac{2}{n}\sum^{\lfloor\sqrt{n}\rfloor}_{l=1}\left(\frac{n}{(l+\alpha)^2}-\frac{n}{(l+1+\alpha)^2}\right)l+\frac{2}{n}\sum^{\lfloor\sqrt{n}\rfloor+1}_{l=2}\left\lbrace\frac{n}{(l+\alpha)^2}\right\rbrace l\\
&-\frac{2}{n}\sum^{\lfloor\sqrt{n}\rfloor}_{l=1}\left\lbrace\frac{n}{(l+1+\alpha)^2}\right\rbrace-\frac{2}{n}\sum^{\lfloor\sqrt{n}\rfloor}_{l=1}\left\lbrace\frac{n}{(l+\alpha)^2}\right\rbrace l\\
&+\frac{2}{n}\left\lbrace \frac{n}{\left(\lfloor\sqrt{n}\rfloor+1+\alpha\right)^2}\right\rbrace(\lfloor\sqrt{n}+1\rfloor)\\
=&\frac{2}{n}\sum^{\lfloor\sqrt{n}\rfloor}_{l=1}\left(\frac{n}{(l+\alpha)^2}-\frac{n}{(l+1+\alpha)^2}\right)l-\frac{2}{n}\left\lbrace\frac{n}{(1+\alpha)^2}\right\rbrace\\
&-\frac{2}{n}\sum^{\lfloor\sqrt{n}\rfloor}_{l=1}\left\lbrace\frac{n}{(l+1+\alpha)^2}\right\rbrace+\frac{2}{n}\left\lbrace \frac{n}{\left(\lfloor\sqrt{n}\rfloor+1+\alpha\right)^2}\right\rbrace(\lfloor\sqrt{n}+1\rfloor)
\end{split}
\end{align*}
As $n$ tends to infinity, the last three terms above become zero, hence we have
\[
L(\alpha) =  \lim_{n \ra +\ns} \frac{2}{n} \sum^{n}_{l=1}\left(\frac{n}{(l+\alpha)^2}-\frac{n}{(l+1+\alpha)^2}\right)l 
=2\sum^{\infty}_{l=1}\left(\frac{1}{(l+\alpha)^2}-\frac{1}{(l+1+\alpha)^2}\right)l.
\]

\hfill $\blacksquare$

\section*{Acknowledgments} The work of PLZ, DH and FJM was supported by a University Staff Doctoral Programme: Building Capacity in Applied Mathematics (USDP-BCAM) grant, ID 32, under the Newton's Operational Development Assistance Fund. The grant is funded by the UK Department for Business, Energy and Industrial Strategy and the Department of Higher Education and Training, South Africa and delivered by the British Council. For further information, please visit https://www.ukri.org/what-we-do/browse-our-areas-of-investment-and-support/newton-fund/.

The research of WMB and RM{\L} was supported by grant no. 2022/47/B/ST1/02114 \emph{Non-random equivalent characterizations of sample boundedness} of National Science Centre, Poland. 

Part of this work was done during RML's stay and PLZ's, DH's and FJM's visit at the African Institute for Mathematical Sciences (AIMS) in Muizenberg, South Africa. AIMS' hospitality is gratefully acknowledged.

\bibliographystyle{alpha} 
\bibliography{/Users/rafallochowski/biblio/biblio}

\begin{thebibliography}{CJMS06}

\bibitem[Bil99]{Billingsley:1999cl}
P. Billingsley.
\newblock {\em Convergence of Probability Measures}.
\newblock John Wiley, New York, 1999.

\bibitem[CJMS06]{Coquet2006}
F. Coquet, A. Jakubowski, J. M{\'e}min, and L.
  S{\l}omi{\'n}ski.
\newblock  Natural Decomposition of Processes and Weak Dirichlet
  Processes, in: {\em In Memoriam Paul-Andr{\'e} Meyer -- S\'eminaire de Probabilit\'es XXXIX}, 81--116,
\newblock Springer Berlin Heidelberg, Berlin, Heidelberg, 2006.

\bibitem[DOS18]{Obloj_local:2015}
M.~Davis, J.~Ob{\l }\'{o}j, and P.~Siorpaes.
\newblock Pathwise stochastic calculus with local times.
\newblock {\em Ann. de l'IHP}, 54(1):1--21, 2018.

\bibitem[F{\"{o}}l81]{Foellmer:1981}
H.~F{\"{o}}llmer.
\newblock Calcul d'{I}t\^{o} sans probabilit\'{e}s.
\newblock {\em S\'{e}minaire de Probabilit\'{e}s XV}, 80:143--150, 1981.

\bibitem[Hir23]{HiraiEJP2023}
Y. Hirai.
\newblock {It\^{o}-F\"{o}llmer calculus in Banach spaces I: the It\^{o}
  formula}.
\newblock {\em Electron. J. Probab.}, 28(none):1--41, 2023.

\bibitem[{\L}OPS21]{LochOblPS:2021}
R.~M. {\L}ochowski, J. Ob{\l}\'oj, D.~J. Pr\"omel, and P.
  Siorpaes.
\newblock Local times and {T}anaka--{M}eyer formulae for c{\`{a}}dl{\`{a}}g
  paths.
\newblock {\em Electron. J. Probab.}, 26(none):1--29, 2021.

\bibitem[Vov15]{Vovk_cadlag:2015}
V.~Vovk.
\newblock {I}t\^{o} calculus without probability in idealized financial
  markets.
\newblock {\em Lith. Math. J.}, 55(2):270--290, 2015.

\bibitem[Zon25]{PhumlaniPhD}
P. L.~Zondi.
\newblock Quadratic variation and stochastic calculus for deterministic paths
  and model-free price paths with jumps.
\newblock {\em PhD Thesis, The University of Limpopo}, 2025.

\end{thebibliography}

\end{document}